\documentclass[10pt]{elsarticle}
\usepackage[cp1251]{inputenc}
\usepackage{amsmath}
\usepackage{amssymb}
\usepackage{amsfonts}
\usepackage{mathtools}
\usepackage{dsfont}
\usepackage{graphicx}
\usepackage{floatflt,epsfig}
\usepackage{lineno,hyperref}
\usepackage{enumerate}
\usepackage{colortbl}
\usepackage{array,tabularx,tabulary,booktabs}
\usepackage{longtable}
\usepackage{multirow}
\usepackage{wrapfig}
\usepackage{subcaption}
\usepackage{caption}
\usepackage[table]{xcolor}
\usepackage{adjustbox}
\usepackage{rotating}
\usepackage{arydshln}
\usepackage{bm}
\usepackage{tabularray}
\newcolumntype{^}{>{\currentrowstyle}}

\usepackage{enumitem}
\usepackage{tabularx}
\usepackage{graphicx}
\usepackage{tikz}
\usepackage{ytableau}
\usepackage{adjustbox}

\usepackage{float}
\journal{Siberian Electronic Mathematical Reports}
\newtheorem{lemma}{Lemma}

\newtheorem{theorem}{Theorem}

\newcommand{\proof}{\medskip\noindent{\bf Proof.~}}
\usepackage{array}
\newcolumntype{P}[1]{>{\raggedright\arraybackslash}p{#1}}

\begin{document}
\renewcommand{\abstractname}{Abstract}
\renewcommand{\refname}{References}
\renewcommand{\arraystretch}{0.9}
\thispagestyle{empty}
\sloppy

\begin{frontmatter}
\title{One more proof about the spectrum of Transposition graph}

\author[01,02]{Artem Kravchuk}
\ead{artemkravchuk13@gmail.com}

\address[01]{Sobolev Institute of Mathematics, Ak. Koptyug av. 4, Novosibirsk 630090, Russia}
\address[02]{Novosibirsk State University, Pirogova str. 2, Novosibirsk, 630090, Russia}

\begin{abstract}
A \emph{Transposition graph} $T_n$ is defined as a Cayley graph over the symmetric group $Sym_n$ generated by all transpositions. This paper shows how the spectrum of $T_n$ can be obtained using the spectral properties of the Jucys-Murphy elements.
\end{abstract}

\begin{keyword}
Transposition graph; integral graph; spectrum; 
\vspace{\baselineskip}
\MSC[2010] 05C25\sep 05E10\sep 05E15
\end{keyword}
% 05C25 graphs and groups
% 05E10 graphs and matrices
% 05E15 combinatorial problems concerning the classical groups;
\end{frontmatter}
%%%%%%%%%%%%%%%%%%%%%%%%%%%%%%%%%%%%%%%%%%%%%%%%%%%%%%%%%%%%%%%%%%%%%%%%%%%%%%%%%%%%%%%%%%%%%%%%%%%%%%%%%%%%%%%%%%%%%%%%%%%%%%%%%%

\section{Introduction}\label{intro}

\emph{Transposition graph} $T_n$ is defined as a Cayley graph over the symmetric group generated by all transpositions. The \emph{spectrum} of a graph is defined as a multiset of distinct eigenvalues of its adjacency matrix together with their multiplicities. The graph $T_n, n\geqslant 2$, is bipartite, so its spectrum is symmetric about zero. The largest eigenvalue of $T_n$ is equal to $\binom{n}{2}$ which implies that all eigenvalues of $T_n$ lie in the interval $[-\binom{n}{2}, \binom{n}{2}]$. In \cite{KY97} it was shown that the eigenvalues of $T_n$ correspond to partitions $\lambda \vdash n$, and formulas for multiplicities of these eigenvalues were given. It was also shown in ~\cite[Lemma 3]{KY97} that $T_n$ is integral which means that its spectrum $Spec(T_n)$ consists of integers. Later and independently in ~\cite[Theorem 2]{KL20}, it was shown that  $T_n$ is integral. The obtaining of the spectrum in \cite{KY97} is based on the Zieschang theorem \cite[Theorem~1]{Z88}. This paper shows how to obtain the spectrum of $T_n$ in an another way, using the spectral properties of the Jucys-Murphy elements. The Jucys-Murphy elements are significant in the study of representation theory related to symmetric groups \cite{OV96} and have diverse applications. For example, spectral properties of the Jucys-Murphy elements have previously been used to find the spectrum of the Star graph~\cite{CF12}. Also in \cite{GKKSV20,GKKSV21}, Jucys-Murphy elements are used to study the spectral properties of this graph.

The paper is organised as follows. First, in Section~\ref{sec2} we give basic facts from the representation theory of the symmetric group that are needed in further proofs. After that in Section~\ref{sec3} we describe the Jucys-Murphy elements and introduce the $J_\Sigma^n$ operator and study it's spectrum. Then in Section~\ref{sec4} we show that the operator $J_\Sigma^n$ is equal to adjacency matrix of $T_n$.

\section{Preliminaries}\label{sec2}
\subsection{Partitions and Young tableaux}

A nonincreasing sequence of positive integers $\lambda = (n_1, \dots, n_k) \vdash n $, for which $n=\sum\limits_{j=1}^kn_j$, is called a \emph{partition} of $n$.

A \emph{Young diagram} is a finite set of cells arranged in left-justified rows, with the lengths of the rows in non-increasing order. It is not hard to see the bijection between Young diagrams and partitions.

\begin{figure}
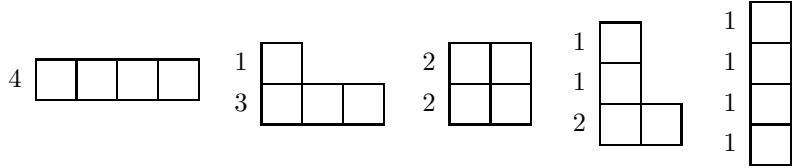

    \centering
\begin{tabular}{@{}m{3cm}@{} @{}m{2.5cm}@{} @{}m{2cm}@{} @{}m{2cm}@{} @{}m{2cm}@{}}
\begin{ytableau}
       \none[4] &  &  & & \\
\end{ytableau}
&
\begin{ytableau}
       \none[1] &   \\
       \none[3] &  &  &  \\
\end{ytableau}
&
\begin{ytableau}
       \none[2] & &  \\
       \none[2] &  &   \\
\end{ytableau}
&
\begin{ytableau}
       \none[1] &  \\
       \none[1] &   \\
       \none[2] &  &  \\
\end{ytableau}
&
\begin{ytableau}
       \none[1] &  \\
       \none[1] &  \\
       \none[1] &    \\
       \none[1] &    \\
\end{ytableau}
\end{tabular}
\caption{Young diagrams corresponding to all partitions of the number 4.}
\label{fig:y_diagr_of_4}
\end{figure}

We say that a Young diagram has the shape $\lambda$, where $\lambda = (n_1, n_2, \dots, n_k)$ if the $i$-th row has $n_i$ cells. We use French notation, in which the rows are constructed from bottom to top. Figure~\ref{fig:y_diagr_of_4} shows Young diagrams corresponding to all partitions of $n=4$.

\begin{figure}
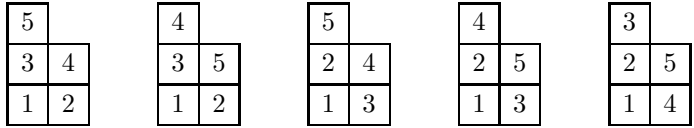

    \centering
\begin{tabular}{@{}m{2cm}@{} @{}m{2cm}@{} @{}m{2cm}@{} @{}m{2cm}@{} @{}m{2cm}@{}}
\begin{ytableau}
       5  \\
       3 & 4 \\
       1 & 2 \\
\end{ytableau}
&
\begin{ytableau}
       4  \\
       3 & 5 \\
       1 & 2 \\
\end{ytableau}
&
\begin{ytableau}
       5  \\
       2 & 4 \\
       1 & 3 \\
\end{ytableau}
&
\begin{ytableau}
       4  \\
       2 & 5 \\
       1 & 3 \\
\end{ytableau}
&
\begin{ytableau}
       3  \\
       2 & 5 \\
       1 & 4 \\
\end{ytableau}
\end{tabular}
\caption{Standard Young tables of the shape $\lambda = (2, 2, 1)$.}
    \label{fig:st_y_t_221}

\end{figure}

A \emph{standard Young tableau} of the shape $\lambda \vdash n$ is the filling of a Young diagram of the shape $\lambda$ with numbers from $1$ to $n$ such that the numbers in each row and each column are in ascending order. Let $\mathcal{T}_\lambda$ be the set of standard Young tableaux of the shape $\lambda\vdash n$. Figure~\ref{fig:st_y_t_221} shows all standard Young tableaux from $\mathcal{T}_{(2, 2, 1)}$.

For each cell at position $(i,j)$, the \emph{hook} $H_{ij}^{\lambda}$ is defined as the set of cells $(k, l)$ such that $k=i$ and $l\geqslant j$ or $l=j$ and $k\geqslant i$. The hook length $h_{ij}^{\lambda}$ is the number of cells in $H_{ij}^{\lambda}$. The number of standard Young tableaux of the shape $\lambda$ is obtained using the hook-length formula~\cite[Theorem 1]{FRT54}:
\begin{equation}\label{hook_formuls}
f_{\lambda} = \frac{n!}{\prod_{(ij)}h_{ij}^{\lambda}}.
\end{equation}

Consider some standard Young tableau $t_{\lambda} \in \mathcal{T_\lambda}, \lambda \vdash n$. For each $i \in \{1, \dots, n\}$ we define $c_{\lambda}^t(i) = y - x$, where $x$ and $y$ are the abscissa and ordinate of number $i$, respectively. For example, for the first table from Figure~\ref{fig:st_y_t_221} we have $c(5) = 3-1 = 1, c(4) = 0, c(3) = 1, c(2) = -1, c(1) = 0.$

Note that the sum $c_{\lambda}^t (i), i \in \{0, \dots, n\}$ does not depend on how the numbers in the standard Young table $t_\lambda$ are filled. Let us denote by $s_{\lambda}$ the following sum: 

\begin{equation}\label{sum_c}
s_{\lambda} = \sum\limits_{i=1}^{n} c_{\lambda}(i).
\end{equation}

For example, for the first Young diagram in Figure~\ref{fig:y_diagr_of_4}, $s_{\lambda}=-6$, and for the Young diagram in Figure~\ref{fig:st_y_t_221}, $s_{\lambda}=3$.

Let $\lambda =(n_1, n_2, \dots, n_k) \vdash n$.
\begin{lemma}\label{sl_expr}
\begin{equation}\label{sl}
    s_{\lambda} = -\sum_{j=1}^{k}\frac{n_j(n_j-2j+1)}{2}.
\end{equation}
\end{lemma}
\proof 

Let us consider the $j$-th row in the Young tableau of the shape $\lambda$. Elements from this line add the expression $\sum\limits_{i=1}^{n_j}(j-i)=-\frac{n_j(n_j-2j+1)}{2}$ to $s_{\lambda}$.
\hfill $\square$

Note that different partitions can correspond to the same value of $s_\lambda$. For example,  $s_{(3, 3)}=s_{(4, 1, 1)} = -3$.

The \emph{conjugate of a partition $\lambda$}, denoted by $\lambda'$, is also a partition of $n$ where its parts are the nonincreasing sequence $\lambda' = (n_1', n_2',\ldots,n_{\kappa}')$, where $n_j' = \sum\limits_{i|n_i \geqslant j}1$.

\begin{lemma}
    $s_{\lambda} =-s_{\lambda'}.$
\end{lemma}

\proof The proof is the same as for Lemma~\ref{sl_expr} if we consider the columns for the Young table corresponding to $\lambda'$.

\subsection{Regular representation}

Let $GL(V)$ represent the collection of all invertible transformations of vector space $V$ onto itself, known as the general linear group of $V.$ A representation of a group $G$ on a vector space $V$ over a field $\mathbb{C}$ is a  homomorphism $\rho: G \longrightarrow GL(V)$. $V$ is also called \emph{$G$-module}.

Let $W$ be a vector subspace of $V$. $W$ is considered $\rho$-invariant if $\rho(g)(w) \in W$ for any $g \in G$ and $w \in W$. A representation $\rho$ is defined as irreducible if the only invariant subspaces of V are $\{0\}$ and $V$. The irreducible representations of $Sym_n$ are indexed by partitions of $n$. We denote by $V_{\lambda}$ the irreducible module associated with the partition $\lambda \vdash n$.

Let $G$ be a finite group of order $n$. Consider a vector space $V$ of dimension $n$, having a basis whose elements correspond to the elements of $G$. We denote this basis as $\{e_g\}_{g\in G}$, where $e_g$ is the basis vector corresponding to the element $g$ of the group $G$. The \emph{regular representation} of the group $G$ is a representation $\rho$, defined as follows:

$$\rho(g)(e_h) = e_{gh}$$
for all $g,h\in G$.

The regular representation of $Sym_n$ can be decomposed into irreducible submodules in the following manner \cite[Prop. 1.10.1]{S01}:

\begin{equation}\label{reg-decomp}
    \mathbb{C}[Sym_n] = \underbrace{\left( V_{\lambda_1} \oplus \cdots \oplus V_{\lambda_1} \right)}_{\dim(V_{\lambda_1})} \oplus \cdots \oplus \underbrace{\left( V_{\lambda_k} \oplus \cdots \oplus V_{\lambda_k} \right)}_{\dim(V_{\lambda_k})} =\bigoplus_{\lambda \vdash n} V_\lambda^{dim(V_{\lambda})}, 
\end{equation}
and $dim(V_\lambda)$ = $f_{\lambda}.$

\section{Jucys-Murphy elements and spectrum of $J_\Sigma^n$}\label{sec3}
The \emph{Jucys-Murphy elements} are elements of group algebra $\mathbb{C}[Sym_n]$. These elements were introduced independently by Jucys in 1966 \cite{J66} and Murphy in 1981 \cite{M81} and are defined as follows:

\begin{equation}
   J_i^n = (1, i) + (2, i) + \dots + (i - 1, i)=\sum\limits_{j=1}^{i}(j, i), i=2, \dots, n, 
\end{equation}
where by the summation of transpositions means the summation of the matrices corresponding to these transpositions in some representation.

The Jucys-Murphy elements can be viewed as a linear operator acting on the basis vectors of the regular representation as follows:

$$J_i^n (e_{\sigma_k}) = \sum\limits_{j=1}^i e_{(j,i)\cdot\sigma_k}.$$

The spectrum of Jucys-Murphy elements is obtained by the following theorem.
\begin{theorem}\label{J-M-th}{\rm \cite{J74}}
Let $\lambda \vdash n$. Then there exists a basis $(v_t)_{t\in \mathcal{T}_\lambda}$ of irreducible module $V_{\lambda}$, indexed by standart Young tables of shape $\lambda$, such that for all $i \in \{2, \dots, n \}$, one has:

\begin{equation}
    J_i^n (v_t) = c^t_{\lambda}(i)v_t.
\end{equation}
\end{theorem}
Let us denote by $J_{\Sigma}^n$ the following expression:

\begin{equation}
    J_{\Sigma}^n = \sum\limits_{i=2}^n J_i^n.
\end{equation}

Consider the basis $(v_t)_{t\in \mathcal{T}_\lambda}$ from Theorem~\ref{J-M-th}. From the definition of $J_{\Sigma}^n$ and $s_\lambda$ it follows directly that $\{s_{\lambda}\}_{\lambda\vdash n}$ are eigenvalues of the operator $J_{\Sigma}^n$:
\begin{equation}
    J_\Sigma^n (v_t) = s_{\lambda}v_t.
\end{equation}

The following lemma gives an expression of the multiplicity of the eigenvalue $s_\lambda$ for the operator $J_\Sigma^n$.

\begin{lemma}
    $mul(s_\lambda) = \sum\limits_{\eta : s_{\alpha} = s_{\lambda}} f_{\alpha}^2$.
\end{lemma}
\proof
By~(\ref{reg-decomp}), the number of irreducible modules $V_\lambda$ is $f_{\lambda}$. In each of $V_{\lambda}$ the number of eigenvectors $v_t$ corresponding to the eigenvalue $s_\lambda$ is equal to $f_\lambda$ too. Note also that $mul(s_{\lambda}) = mul(s_{\lambda'})$.
\hfill $\square$

% Thus, we can describe the spectrum of the operator $J_{\Sigma}^n$ as follows. To each partition $\lambda = (n_1, \dots, n_k)$ corresponds an eigenvalue $s_\lambda =  -\sum_{j=1}^{k}\frac{n_j(n_j-2j+1)}{2}$ with multiplicity $\sum\limits_{\eta : s_{\eta} = s_{\lambda}} f_{\eta}^2$. 

We denote by $\eta_{\lambda} = s_{\lambda'} = -s_{\lambda}$. Then we are able to formulate the following theorem. 

\begin{theorem}\label{js-th}
    The spectrum of  $J_\Sigma^n$ consists of the eigenvalues corresponding to the partitions $\lambda\vdash n$. Any partition $\lambda \vdash n$ corresponds to an eigenvalue $\eta_\lambda$ of  $J_\Sigma^n$ by the following expression:
    \begin{equation}\label{eig-expr}
        \eta_{\lambda} = \sum_{j=1}^{k}\frac{n_j(n_j-2j+1)}{2}
    \end{equation}
    with multiplicity $mul(\eta_\lambda) =  \sum\limits_{\alpha : \eta_{\alpha} = \eta_{\lambda}} f_{\alpha}^2$.
\end{theorem}

\section{Equivalence of $J_\Sigma^n$ and the adjacency matrix of  $T_n$}\label{sec4}
% Note that $J_\Sigma^n$ is the sum of all transpositions:

% \begin{equation}
%     J_\Sigma=\sum\limits_{i<j}(i, j).
% \end{equation}
Let us denote by $A_{T_n}$ the adjacency matrix (or equal linear operator) of $T_n$. Then we can formulate the following lemma.

\begin{lemma}\label{eq-lemma}
   $A_{T_n} = J_{\Sigma}^n.$ 
\end{lemma}
\proof
Let us see how the operator $A_{T_n}$ acts on an arbitrary basis element $e_{\sigma_k}$:

\begin{equation}\label{adj_1}
    A_{T_n} (e_{\sigma_k}) = \sum\limits_{\sigma_l: (\sigma_l, \sigma_k) \in E(T_n)} e_{\sigma_l}.
\end{equation}
An edge $(\sigma_l, \sigma_k) \in E(T_n)$ if and only if $\sigma_l = (i, j)\sigma_k$ for any transposition $(i, j)$. Therefore,~(\ref{adj_1}) can be rewritten in the following form:
\begin{equation}
    A_{T_n}(e_{\sigma_k}) = \sum\limits_{i<j}e_{(j,i)\sigma_k}.
\end{equation}
Now let us consider the action of $J_\Sigma^n$ on $e_{\sigma_k}$:
\begin{equation}
    J_\Sigma^n(e_{\sigma_k}) = \sum\limits_{i=2}^n J_i^n(e_{\sigma_k}) = \sum_{i=2}^{n}\sum\limits_{j=1}^i e_{(j,i)\sigma_k} = \sum\limits_{j<i}e_{(j,i)\sigma_k}.
\end{equation}
Thus, the action of operators $J_\Sigma^n$ and $A_{T_n}$ on an arbitrary basis element $e_{\sigma_k}$ is the same, so operators $J_{\Sigma}^n$ and $A_{T_n}$ are equal to each other.
\hfill $\square$

Lemma~\ref{eq-lemma}  shows that the spectrum of $T_n$ can be expressed using Theorem~\ref{js-th}.

\section{Discussions and further research}

In this paper, it was shown how one can use the spectral properties of Jucys-Murphy elements to obtain the spectrum of $T_n$.

However, by looking at the expression (~\ref{eig-expr}) it is not possible to answer the following questions:

\begin{itemize}
    \item  are all integers from the interval $[0, k], k\in \mathbb{N}$, lie in the spectrum of $T_n$?
    \item what is the asymptotics for the number of unique values in the spectrum of $T_n$?
\end{itemize}

Ideally, one would like to be able to answer the question whether a given integer number $k \in [-\binom{n}{2}, \binom{n}{2}]$ lies in the spectrum of $T_n$. The papers \cite{KK22,KK23, K24} partially address these issues.

\section*{Acknowledgements}
 The author was supported by the Mathematical Center in Akademgorodok, under agreement No. 075-15-2022-281 with the Ministry of Science and High Education of the Russian Federation.

\end{document}